\documentclass[12pt]{amsart}
\usepackage{cases}
\usepackage{txfonts}
\usepackage{latexsym}
\usepackage{amsthm}
\usepackage{mathrsfs}
\usepackage{amssymb, amsmath}

\def\opn#1#2{\def#1{\operatorname{#2}}} 
\opn\chara{char} \opn\length{\ell}
\opn\projdim{proj\,dim} \opn\injdim{inj\,dim} \opn\rank{rank}
\opn\depth{depth} \opn\grade{grade} \opn\height{height}
\opn\embdim{emb\,dim} \opn\codim{codim}

\opn\Tr{Tr} \opn\bigrank{big\,rank}
\opn\superheight{superheight}\opn\lcm{lcm}
\opn\trdeg{tr\,deg}%
\opn\reg{reg} \opn\lreg{lreg}
\opn\div{div} \opn\Div{Div} \opn\cl{cl} \opn\Cl{Cl}
\opn\Spec{Spec} \opn\Supp{Supp} \opn\supp{supp} \opn\Sing{Sing}
\opn\Ass{Ass}
\opn\Ann{Ann} \opn\Rad{Rad} \opn\Soc{Soc}
\opn\Ker{Ker} \opn\Coker{Coker} \opn\Im{Im} \opn\Hom{Hom}
\opn\Tor{Tor} \opn\Ext{Ext} \opn\End{End} \opn\Aut{Aut} \opn\id{id}

\opn\nat{nat}
\opn\pff{pf}
\opn\Pf{Pf} \opn\GL{GL} \opn\SL{SL} \opn\mod{mod} \opn\ord{ord}
\opn\aff{aff} \opn\con{conv} \opn\relint{relint} \opn\st{st}
\opn\lk{lk} \opn\cn{cn} \opn\core{core} \opn\vol{vol}
\opn\gr{gr}

\def\pot#1#2{#1[\kern-0.28ex[#2]\kern-0.28ex]}

\opn\dirlim{\underrightarrow{\lim}}
\opn\invlim{\underleftarrow{\lim}}
\newtheorem{theorem}{Theorem}

\newtheorem{lemma}[theorem]{Lemma}

\def\qed{\ifhmode\textqed\fi
   \ifmmode\ifinner\quad\qedsymbol\else\dispqed\fi\fi}
\def\textqed{\unskip\nobreak\penalty50
    \hskip2em\hbox{}\nobreak\hfil\qedsymbol
    \parfillskip=0pt \finalhyphendemerits=0}
\def\dispqed{\rlap{\qquad\qedsymbol}}

\opn\ini{in} \opn\inm{inm} \opn\Sym{Sym}

\begin{document}
\title{distinct zeros of the riemann zeta-function}
\author{XiaoSheng Wu}
\date{}
\address {School of Mathematical Sciences\\ University of Science and Technology of China\\ Hefei 230026, P. R.
China.}
\email {xswu12@ustc.edu.cn}

\begin{abstract} In this paper, we prove that there are more than 66.036\% of zeros of the Riemann zeta-function are distinct.
\end{abstract}
\maketitle


\section{Introduction}
Let $\zeta(s)$ be the Riemann zeta-function, where $s=\sigma+it$. It is defined for $\sigma>1$ by
\begin{align}
  \zeta(s)=\sum_{n\ge1}n^{-s}.\notag
\end{align}
The Riemann-von Mangoldt formula \cite{Tic} states that $N(T)$, the number of non-trivial zeros $\rho=\beta+i\gamma$ of $\zeta(s)$ with $0<\gamma\le T$, satisfies
\begin{align}
   &N(T)=\frac{T}{2\pi}\log\frac{T}{2\pi e}+\frac78+S(T)+O(\frac1T),\notag\\
   &S(T)=\frac1\pi\arg\zeta(\frac12+iT)\ll\log T.\notag
\end{align}
It is generally believed that all the zeros of $\zeta(s)$ are simple (also distinct), which is known as the Simple Zero Conjecture. We define the number of simple zeros and the number of distinct zeros for the Riemann zeta-function as follows
\begin{align}
   &N_d(T)=|\{\rho=\beta+i\gamma: 0<\gamma\le T, \zeta(\rho)=0\}|,\notag\\
   &N_s(T)=|\{\rho=\beta+i\gamma: 0<\gamma\le T, \zeta(\rho)=0, \zeta'(\rho)\neq0\}|.\notag
\end{align}
The Simple Zero Conjecture means $N_d(T)=N_s(T)=N(T)$.

Due to Levinson's method, it is known that more than two-fifths of the zeros are simple (see \cite{BCY, Con-More}). In 1995, Farmer \cite{Far} introduced a combination method, which is based on proportions of simple zeros of $\xi^{(n)}(s)$, $n\ge0$, and proved that at least 63.9\% of the zeros of the Riemann zeta-function are distinct.

In this paper, by the method introduced in our work \cite{Wu}, we prove that there are more than 66.03\% of zeros of the Riemann zeta-function are distinct.
\begin{theorem}
\label{thm1}
For $T$ sufficiently large, we have
\begin{align}
   N_d(T)\ge0.66036N(T).\notag
\end{align}
\end{theorem}

In this paper we always assume that $T$ is a large parameter and $L=\log T$. The number of additional zeros of a analytic function $f$ caused by multiplicity means the number of zeros of $f$ counted according to multiplicity minus one.

We note that if $\rho$ is a non-simple zero of $\zeta(s)$, it must be a zero of
\begin{align}
\label{2.1}
   G(s)=\zeta(s)\psi_1(s)+\zeta'(s)\psi_2(s)
\end{align}
with multiplicity reduced by at most one, here $\psi_1(s)$ and $\psi_2(s)$ can be any analytic function. Thus  the number of additional zeros of $\zeta(s)$ caused by multiplicity in any region is not more than the number of zeros of $G(s)$ in the same region. We partition the whole plane into the left part ($\text{Re}(s)<1/2$) and the right part ($\text{Re}(s)\ge1/2$). Then we will evaluate the number of additional zeros of $\zeta(s)$ caused by multiplicity in each side with different $G(s)$.

To the left side, we choose $G(s)=\xi'(s)$, where
\begin{align}
   \xi(s)=H(s)\zeta(s)\notag
\end{align}
with
\begin{align}
\label{+2}
   H(s)=1/2s(s-1)\pi^{-s/2}\Gamma(s/2).
\end{align}
It is known that for any integer $n\ge0$ the number of zeros for $\xi^{(n)}(s)$ with $0<t<T$ is $N(T)+O(\log T)$ (to see \cite{Con-Zero, Lev-Zero}).
The functional equation for $\zeta(s)$ says that
\begin{align}
\label{+1}
   \xi(s)=\xi(1-s).
\end{align}
Differentiating both side of the above formula in $s$ we have
\begin{align}
   \xi'(s)=-\xi'(1-s).\notag
\end{align}
Hence if $\rho$ is a zero of $\xi'(s)$, so is $1-\rho$.  Let $N_{\xi',c}(T)$ be the number of zeros of $\xi'(1/2+it)$ with $0<t<T$. Then it follows from the symmetry of zeros of $\xi'(s)$ that the number of additional zeros of $\zeta(s)$ in the left side is not more than
\begin{align}
   \frac12\Big(N(T)-N_{\xi',c}(T)\Big)+O(\log T).\notag
\end{align}

To the right side, we choose $G(s)$ as in (\ref{2.1}) with
\begin{align}
\label{+3}
   &\psi_1(s)=\sum_{n\le y}\frac{\mu(n)}{n^{s+R/L}}P_1(\frac{\log y/n}{\log y}),\notag\\
   &\psi_2(s)=\frac1L\sum_{n\le y}\frac{\mu(n)}{n^{s+R/L}}P_2(\frac{\log y/n}{\log y}),
\end{align}
where $\mu$ is the Mobius function, $y=T^\theta$ with $0<\theta<4/7$. Here $P_1,~P_2$ are polynomials with $P_1(0)=P_2(0)$, $P_1(1)=1$ which will be specified later (to see \cite{BCY, Feng, Lev-Dedu} for any more choice of $\psi_i$).

Let $D$ be the closed rectangle with vertices $1/2+it_0$, $3+it_0$, $1/2+iT$, $3+iT$. Here $t_0\le2$ is a given positive constant that is not a ordinate of any zero of $G(s)$ in the region $0<\text{Re}(s)<3,~\text{Im}(s)>0$. Let $N_G(D)$ denote the zeros of $G(s)$ in $D$, including zeros on the left boundary.

Since the number of additional zeros of $\zeta(s)$ caused by multiplicity is not more than the number of zeros of $\xi'(s)$ in the left side and not more than $N_G(D)$ in the right side, we may have the following formula about the number of distinct zeros of $\zeta(s)$
\begin{align}
\label{2.2}
   N_d(T)\ge \frac12N(T)+\frac12N_{\xi',c}(T)-N_G(D).
\end{align}
It is therefore important to give an upper bound for $N_G(D)$ and a lower bound for $N_{\xi',c}(T)$. We will obtain a an upper bound for $N_G(D)$ in section \ref{sec2} and a lower bound better than the known result now for $N_{\xi',c}(T)$ in section \ref{sec3}.

\section{upper bound for $N_G(D)$}
\label{sec2}
An upper bound for $N_G(D)$ can be found in a familiar way by applying  Littlewood's formula (to see \S9.9 of \cite{Tic}). Let \begin{align}
   \sigma_0=1/2-R/L,\notag
\end{align}
where $R$ is a constant to be specified precisely later. Let $D_1$ be the closed rectangle with vertices $\sigma_0+it_0$, $3+it_0$, $\sigma_0+iT$, $3+iT$. Suppose $G(3+it)\neq0$. Determine $\text{arg}G(\sigma+iT)$ by continuation left from $3+iT$ and $\text{arg}G(\sigma+it_0)$ by  continuation left from $3+iTt_0$. If a zero is reached on the upper edge, use $\lim G(\sigma+iT+i\epsilon)$ as $\epsilon\rightarrow+0$. Make horizontal cuts in $D_1$ from the left side to the zeros of $G$ in $D_1$. Applying the Littlewood's formula, we have
\begin{align}
\label{2.4}
   &\int_{t_0}^T\log|G(\sigma_0+it)|dt-\int_{t_0}^T\log|G(3+it)|dt\notag\\
   &\ \ \ \ \ \ \ \ \ \ \ \ \ \ \ \ +\int_{\sigma_0}^3\text{arg}G(\sigma+iT)d\sigma-\int_{\sigma_0}^3\text{arg}G(\sigma+it_0)d\sigma\notag\\
   &\ \ \ \ \ \ \ =2\pi\sum_{\rho\in D_1}\text{dist}(\rho),
\end{align}
where $\text{dist}(\rho)$ is the distance of $\rho$ from the left side of $D_1$.

Recall the definition of $\psi_i$ for $i=1,2$. A direct calculation shows that $\psi_i(s)\ll T$ for $\text{Re}(s)>0$. Hence $G(s)\ll T^2$ for $\text{Re}(s)>0$. Then using Jensen's theorem in a familiar way as in \S9.4 of \cite{Tic}, we have
\begin{align}
\label{2.4.1}
   \int_{\sigma_0}^3\text{arg}G(\sigma+iT)d\sigma=O(\log T).
\end{align}
For $t_0$ is a given constant, it is easy to see
\begin{align}
\label{2.4.2}
   \int_{\sigma_0}^3\text{arg}G(\sigma+it_0)d\sigma=O(1).
\end{align}
By a direct calculation we can see
\begin{align}
   \zeta'(3+it)\psi_i(3+it)\ll O(1/\log T).\notag
\end{align}
Hence we have from (\ref{2.1}) that
\begin{align}
\label{2.5}
   \int_{t_0}^T\log |G(3+it)|dt=\int_{t_0}^T\log |\zeta(3+it)\psi_i(3+it)|dt+O(T/\log T).
\end{align}
Since for $\sigma>1$
\begin{align}
   \log\zeta(s)=-\sum_{n\ge1}\frac{\Lambda(n)}{n^s\log n},\notag
\end{align}
it follows taking the real part that
\begin{align}
\label{2.6}
   \int_{t_0}^T\log |\zeta(3+it)|dt\ll1.
\end{align}
For the entire function $\psi_i(s)$, it is easy to see, for $\sigma\ge3$,
\begin{align}
   |\psi_i(s)-1|\le\frac{1}{2^\sigma}+\frac1{3^\sigma}+\int_3^\infty\frac{\nu}{\nu^\sigma}\le\frac{1}{2^\sigma}+\frac52\frac{1}{3^\sigma}<2^{1-\sigma}.\notag
\end{align}
Therefore, $\log\psi_i(s)$ is analytic for $\sigma\ge3$. Integrating on the contour $\sigma+iT,~3\le\sigma<\infty$; $3+it,~t_0\le t\le T;~\sigma+iT,~3\le \sigma<\infty$ gives
\begin{align}
\label{2.7}
   \int_{t_0}^T\log|\psi_1(3+it)|dt\le\bigg|\int_{t_0}^T\log\psi_1(3+it)dt\bigg|\le8\int_3^\infty\frac{d\sigma}{2^\sigma}=O(1).
\end{align}
Substituting (\ref{2.6}), (\ref{2.7}) into (\ref{2.5}) we have
\begin{align}
\label{2.4.3}
   \int_{t_0}^T\log |G(3+it)|dt\ll T/\log T.
\end{align}
Then using (\ref{2.4.1}), (\ref{2.4.2}), (\ref{2.4.3}) in (\ref{2.4}), we have
\begin{align}
   \int_{t_0}^T\log|G(\sigma_0+it)|dt+O(T/\log T)=2\pi\sum_{\rho\in D}\text{dist}(\rho).\notag
\end{align}
Since all the zeros of $G$ in closed rectangle $D$ are at least distance $1/2-\sigma_0$ from $\sigma=\sigma_0$, it follows that
\begin{align}
\label{+4}
  2\pi(1/2-\sigma_0)N_G(D)\le\int_{t_0}^T\log|G(\sigma_0+it)|dt+O(T/\log T).
\end{align}
Using the concavity of the logarithm,
\begin{align}
   \int_{t_0}^T\log|G(\sigma_0+it)|dt&=1/2\int_{t_0}^T\log|G(\sigma_0+it)|^2dt\notag\\
   &\le\frac12T\log\Big(\frac1T\int_{t_0}^T|G(\sigma_0+it)|^2dt\Big).\notag
\end{align}
Substituting $\sigma_0=1/2-R/L$, then we have from (\ref{+4}) that
\begin{align}
\label{2.11}
   N_G(D)\le\frac{TL}{4\pi R}\log\Big(\frac1T\int_{t_0}^T|G(\sigma_0+it)|^2dt\Big).
\end{align}
Hence an upper bound for $N_G(D)$ may be obtained by evaluating the mean value integral of $G$ on the $\sigma_0$-line.

To evaluate the mean value integral of $G$ on the $\sigma_0$-line, we need the following two Lemmas.
\begin{lemma}
\label{lem1}
Suppose that $\delta>0$ and $\Delta=T^{1-\delta}$. Then
\begin{align}
   \frac1T\int_{t_0}^T|G(\sigma_0+it)|^2dt=\frac{1}{\Delta\pi^{1/2}}\int_{-\infty}^{\infty}e^{-(t-w)^2\Delta^{-2}}|G(\sigma_0+it)|^2dt+o_\delta(1)\notag
\end{align}
uniformly for $T\le w\le2T$.
\end{lemma}
This lemma follows exactly as in Section 3 of Balasubramanian, Conrey, and Heath-Brown \cite{BCH}.

\begin{lemma}
\label{lem2}
Let $a,b\in \mathbb{C}$ with $a,b\ll1$, and put $\alpha=a/L,~\beta=b/L$ where $L=\log T$. Let $s_0=1/2+iw$ with $T\le w\le2T$. Suppose that $\delta>0,~\Delta=T^{1-\delta}$ and that $y=T^\theta$ with $0<\theta<4/7$. For $i,j=1,2$, let
\begin{align}
   &g(a,b,w,P_1,P_2)\notag\\
   &\ \ \ \ \ \ =\frac{1}{i\Delta\pi^{1/2}}\mathop{\int}_{(1/2)}e^{(s-s_0)^2\Delta^{-2}}\zeta(s+\alpha)\zeta(1-s+\beta)\psi_1(s-R/L)\psi_2(1-s-R/L)ds,\notag
\end{align}
where $(1/2)$ denotes the straight line path from $1/2-i\infty$ to $1/2+i\infty$. Then
\begin{align}
\label{3.1}
   g(a,b,w,P_i,P_j)=\frac{\Sigma(b,a,P_i,P_j)-e^{-a-b}\Sigma(-a,-b,P_i,P_j)}{\theta(a+b)}+o_\delta(1).
\end{align}
uniformly in $a,~b$, and $w$, where
\begin{align}
\label{3.2}
   \Sigma(a,b,P_i,P_j)=\int_0^1(P_i'(x)+a\theta P_i(x))(P_j'(x)+b\theta P_j(x))dx.
\end{align}
\end{lemma}
This lemma is the Lemma 2 of Conrey \cite{Con-More}.

We now evaluate the mean value integral of $G$ on the $\sigma_0$-line . From Lemma \ref{lem1} we have
\begin{align}
     \frac1T\int_{t_0}^T|G(\sigma_0+it)|^2dt&=\frac{1}{\Delta\pi^{1/2}}\int_{-\infty}^{\infty}e^{-(t-w)^2\Delta^{-2}}|G(\sigma_0+it)|^2dt+O_\delta(1)\notag\\
     &=\frac{1}{i\Delta\pi^{1/2}}\mathop{\int}_{(1/2)}e^{(s-s_0)^2\Delta^{-2}}G(s-R/L)G(1-s-R/L)ds+o_\delta(1).\notag
\end{align}
Then recalling the definition of $G$ in (\ref{2.1}) we have
\begin{align}
     \frac1T\int_{t_0}^T|G(\sigma_0+it)|^2dt=&\frac{1}{i\Delta\pi^{1/2}}\mathop{\int}_{(1/2)}e^{(s-s_0)^2\Delta^{-2}}\zeta\psi_1(s-R/L)\zeta\psi_1(1-s-R/L)ds\notag\\
     +&\frac{1}{i\Delta\pi^{1/2}}\mathop{\int}_{(1/2)}e^{(s-s_0)^2\Delta^{-2}}\zeta\psi_1(s-R/L)\zeta'\psi_2(1-s-R/L)ds\notag\\
     +&\frac{1}{i\Delta\pi^{1/2}}\mathop{\int}_{(1/2)}e^{(s-s_0)^2\Delta^{-2}}\zeta'\psi_2(s-R/L)\zeta\psi_1(1-s-R/L)ds\notag\\
     +&\frac{1}{i\Delta\pi^{1/2}}\mathop{\int}_{(1/2)}e^{(s-s_0)^2\Delta^{-2}}\zeta'\psi_2(s-R/L)\zeta'\psi_2(1-s-R/L)ds+o_\delta(1).\notag
\end{align}
We may evaluate every item in the above formula by using Lemma \ref{lem2}. Then
\begin{align}
\label{3.3}
     \frac1T\int_{t_0}^T|G(\sigma_0+it)|^2dt&=g(a,b,w,P_1,P_1)\big|_{a=b=-R}+\partial_bg(a,b,w,P_1,P_2)\big|_{a=b=-R}\notag\\
     &+\partial_ag(a,b,w,P_2,P_1)\big|_{a=b=-R}+\partial_a\partial_bg(a,b,w,P_2,P_2)\big|_{a=b=-R}+o_\delta(1).
\end{align}
Taking $\theta=4/7-\epsilon$, $R=1.023$,
\begin{align}
   &P_1(x)=x-0.064x(1-x)+0.112x^2(1-x),\notag\\
   &P_2(x)=1.305x-0.276x^2-0.025x^3,\notag
\end{align}
and making $\epsilon\rightarrow0$ in (\ref{3.3}), we get from (\ref{2.11}) that
\begin{align}
   N_G(D)\le0.27442N(T).\notag
\end{align}

\section{zeros of $\xi'(s)$ on the critical line (1/2-line)}
\label{sec3}
The fact that a positive proportion of zeros for $\xi'(s)$ are on the critical line was first proved by Levinson. Using the method introduced in his work \cite{Lev-More}, Levinson \cite{Lev-Zero} proved that at least 71\% zeros of $\xi'(s)$ are on the critical line.

In 1983, Conrey \cite{Con-Zero} made a careful study of zeros of $\xi^{(n)}(s)$ on the critical line and proved that at least 81.37\% zeros of $\xi'(s)$ are on the critical line, and later, in his work \cite{Con-More} to prove more than two fifths of the zeros of the Riemann zeta-function are on the critical line, he successfully proved that the mollifier with length $\theta=4/7-\epsilon$ is available. It is obvious that this result can be used in \cite{Con-Zero}. Using the mollifier with length $\theta=4/7-\epsilon$ in \cite{Con-Zero}, one may obtain that at least 82.402\% zeros of $\xi'(s)$ are on the critical line.

We note that if not following the way in \cite{Con-Zero} but using the Levinson's method generalized by Conrey in \cite{Con-More}, we may obtain a better result on this problem. From the functional equation (\ref{+1}) it is easy to see that $\xi^{(n)}(s)$ is real for $s=1/2+it$ when $n$ is even and is purely imaginary when $n$ is odd. Let $\delta\neq0$ be real, $g_n, ~n\ge1$, be complex numbers with $g_n$ real if $n$ is even and $g_n$ purely imaginary if $n$ is odd.
Now define
\begin{align}
   \eta(s)=(1-\delta)\xi(s)+\delta\xi'(s)L^{-1}+\sum_{n=1}^Ng_n\xi^{(n)}(s)L^{-n}\notag
\end{align}
for some fixed $N$. Then, for $s=1/2+it$,
\begin{align}
   \delta \xi'(s)=\text{Im}\eta(s),\notag
\end{align}
so that $\xi'(s)=0$ on $\sigma=1/2$ if and only if $\text{Im}\eta(s)=0$. Observe that for every change of $\pi$ in the argument of $\eta(s)$ it must be the case that $\text{Im}\eta(s)$ has at least one zero. Hence it follows that
\begin{align}
\label{s3.3}
   N_{\xi',c}(T)\ge\frac1\pi\Delta_C\text{arg}\eta(s),
\end{align}
where $\Delta_C\text{arg}$ stands for the variation of the argument as $s$ runs over the critical line from $1/2+it_0$ to $1/2+iT$ passing the zeros of $\eta(s)$ from the east side.

To estimate the change in argument of $\eta(s)$ on the critical line, we let $\eta(s)=H(s)V(s)$, where $H(s)$ is defined in (\ref{+2}) and
\begin{align}
   V(s)=(1-\delta)\zeta(s)+\frac{\delta}{L}\bigg(\frac{H'}{H}(s)\zeta(s)+ \zeta'(s)\bigg)+\sum_{n=1}^N\frac{g_n}{L^{n}}\sum_{k=0}^n\binom{n}{k}\frac{H^{(n-k)}(s)}{H(s)}\zeta^{(k)}(s). \notag
\end{align}
By the Stirling formula, for $|t|\ge2$, we have
\begin{align}
   \text{arg}H(1/2+it)=\frac t2\log\frac{|t|}{2\pi e}+O(1),\notag
\end{align}
and
\begin{align}
   \frac{H^{(m)}}{H}(s)=\Big(\frac12\log\frac {s}{2\pi}\Big)^m(1+O(1/|t|))\notag
\end{align}
for $t\ge10$, $0<\sigma<A_1$, here $A_1$ can be any positive constant. (For a proof of these formulas, see Lemma 1 of \cite{Con-Zero}.) Hence we may have
\begin{align}
\label{s3.2}
   \Delta\text{arg}\eta(1/2+it)\big|_{t_0}^T=T\log\frac{T}{2\pi e}+\Delta\text{arg}V(1/2+it)\big|_{t_0}^T+O(T)
\end{align}
and denote $V(s)$ by
\begin{align}
   V(s)=\Bigg\{\Bigg(1-\delta+\delta\bigg(\frac{\log \frac{s}{2\pi}}{2L}+\frac1{L}\frac{d}{ds}\bigg)Q_0\bigg(\frac{\log \frac{s}{2\pi}}{2L}+\frac1{L}\frac{d}{ds}\bigg)\Bigg)\zeta(s)\Bigg\}(1+O(1/|t|))\notag
\end{align}
with
\begin{align}
   Q_0(x)=1+\sum_{n=1}^N\frac{g_n}{\delta}x^{n-1}.\notag
\end{align}

As in (\ref{+3}), we use the mollifier
\begin{align}
\label{s+1}
   \psi(s)=\sum_{n\le y}\frac{\mu(n)}{n^{s+R/L}}P\bigg(\frac{\log y/n}{\log y}\bigg),
\end{align}
where $\mu$ is the Mobius function, $y=T^\theta$ with $0<\theta<4/7$. Here $P$ is a polynomial with $P(0)=0$, $P(1)=1$ which will be specified later.
By the Cauchy's argument principle, it is not difficult to see
\begin{align}
    \Big|\Delta\text{arg}V(1/2+it)\big|_{t_0}^T\Big|=2\pi N_{\xi'}(D)(1+o(1)).\notag
\end{align}
Here $N_{\xi'}(D)$ denote the number of $\xi'(s)$ in the closed rectangle $D$ defined before.
Then, if $Q_0(1/2)=2$, by applying Jenson's theorem and Littlewood's formula as in section \ref{sec2}, we can show that
\begin{align}
\label{s3.1}
   |\Delta\text{arg}V(1/2+it)|_{t_0}^T|\le \frac{L}{R}\int_{t_0}^T\log|V\psi(\sigma_0+it)|dt(1+o(1)),
\end{align}
where
\begin{align}
   \sigma_0=1/2-R/L,\notag
\end{align}
$R\ll1$ is a positive real number. Here, the reason for requiring the condition $Q_0(1/2)=2$ is to ensure the integration $\int_{t_0}^T\log |V\psi(3+it)|dt$ caused by using Littlewood's formula is $O(T/L)$. To evaluate the integral in the right of (\ref{s3.1}) we use the following useful approximation to $V(s)$,
\begin{align}
\label{s3.4}
   U(s)=\bigg(1-\delta+\delta\Big(1+\frac2{L}\frac{d}{ds}\Big)Q\Big(-\frac1{L}\frac{d}{ds}\Big)\bigg)\zeta(s),
\end{align}
where
\begin{align}
   Q(x)=\frac12Q_0(1/2-x).\notag
\end{align}
If restrict $Q(x)$ be real polynomial, then the restriction of $g_n$ and the condition that $Q_0(1/2)=2$ are equivalent to $Q'(x)=Q'(1-x)$ and $Q(0)=1$. It is easy to see that the error caused by the substitution of $V$ with $U$ can be absorbed by the error term in (\ref{s3.1}). Hence we have from (\ref{s3.3})-(\ref{s3.4})that
\begin{align}
   N_{\xi',c}(T)\ge N(T)-\frac{L}{\pi R}\int_{t_0}^T\log|U\psi(\sigma_0+it)|dt(1+o(1))+O(L).\notag
\end{align}
Recall that $\sigma_0=1/2-R/L$. From the concavity of the logarithm we have
\begin{align}
\label{s3.4}
   N_{\xi',c}(T)\ge N(T)-\frac{TL}{2\pi R}\log\Big(\frac1T\int_{t_0}^T|U\psi(\sigma_0+it)|^2dt\Big).
\end{align}
By Lemma \ref{lem1} and Lemma \ref{lem2} we have
\begin{align}
\label{s3.11}
   \frac1T\int_{t_0}^T|U\psi(\sigma_0+it)|^2\sim&\bigg(1-\delta+\delta(1+2\partial_a)Q(-\partial_a)\bigg)\bigg(1-\delta+\delta(1+2\partial_b)Q(-\partial_b)\bigg)\notag\\
   &\times\Bigg(\frac{g(b,a)-e^{-a-b}g(-a,-b)}{\theta(a+b)}\Bigg)\Bigg|_{a=b=-R},
\end{align}
where
\begin{align}
   g(a,b)=\int_0^1\Big(P'(t)+a\theta P(t)\Big)\Big(P'(t)+b\theta P(t)\Big)dt.\notag
\end{align}

Let $\theta=4/7-\epsilon,~R=1.104,~\delta=0.869$. Taking
\begin{align}
   P(x)= x-0.274x(1-x)-0.334x^2(1-x)+0.005x^3(1-x),\notag
\end{align}
\begin{align}
   Q(x)=1-0.609x-0.572(x^2/2-x^3/3)-4.895(x^3/3-x^4/2+x^5/5)\notag
\end{align}
into (\ref{s3.11}) and making $\epsilon\rightarrow0$, then we have by (\ref{s3.4}) that
\begin{align}
    N_{\xi',c}(T)\ge0.86957 N(T).\notag
\end{align}

It is easy to see that the way in this section can also give better numerical results about the proportion of zeros of $\xi^{(n)}(s)$, $n\ge2$ on the critical line, however one may find that this way is useless when consider simple zeros of $\xi^{(n)}(s)$, $n\ge1$ on the critical line.

\section{completion of the proof}
We have obtained that
\begin{align}
  N_G(D)\le 0.27442N(T)\notag
\end{align}
in section \ref{sec2} and
\begin{align}
    N_{\xi',c}(T)\ge0.86957 N(T).\notag
\end{align}
in section \ref{sec3}.
Then by (\ref{2.2}) we have
\begin{align}
   N_d(T)&\ge \frac12N(T)+\frac12N_{\xi',c}(T)-N_G(D)\notag\\
   &\ge(\frac12+0.434785-0.27442)N(T)>0.66036N(T),\notag
\end{align}
Hence we have proved Theorem \ref{thm1}.

\medskip
\medskip
\noindent

\medskip
\noindent


\begin{thebibliography}{99}

\bibitem{BCH} R. Balasubramanian, J. B. Conrey, and D. R. Heath-Brown, `Asymptotic mean square of the product of the Riemann zeta-function and a Dirichlet Polynomial,' \textit{J. Reine Angew. Math.,} 357 (1985), 161-181.
\bibitem{BCY} H. Bui, J. B. Conrey, and M. Young, `More than 41\% of the zeros of the zeta function are on the critical line,' \textit{Acta Arith.,} 150 (2011), no.1, 35-64.
\bibitem{Con-Zero} J. B. Conrey, `Zeros of derivatives of the Riemann's $\xi$-function on the critical line,' \textit{J. Number Theory,} 16 (1983), 49-74.
\bibitem{Con-More} J. B. Conrey, `More than two fifths of the zeros of the Riemann zeta function are on the critical line,' \textit{J. Reine Angew. Math.,} 339 (1989), 1-26.
\bibitem{Feng} S. Feng, `Zeros of the Riemann zeta function on the critical line,' \textit{J. Number Theory,} 132 (4) (2012), 511-542.
\bibitem{Far} D. Farmer, `Counting distinct zeros of the Riemann zeta-function,' \textit{Electron. J. Combin.,} 2 (1995), Research Paper 1, approx. 5 pp.
\bibitem{Lev-More} N. Levinson, `More than one third of zeros of Riemann's zeta-function are on $\sigma=1/2$,' \textit{Adv. Math.,} 13 (1974), 383-436.
\bibitem{Lev-Zero} N. Levinson, `Zeros of derivative of Riemann's $\xi$-function,' \textit{Bull. Amer. Math. Soc.,} 80, No.5 (1974), 951-954.
\bibitem{Lev-Dedu} N. Levinson, `Deduction of semi-optimal mollifier for obtaining lower bounds for $N_0(T)$ for Riemann's zeta function,' \textit{Proc. Natl. Acad. Sci. USA,} 72 (1975), 294-297.
\bibitem{Tic} E. C. Tichmarsh, The theory of the Riemann zeta-function.
          revised by D. R. HeathBrown, Clarendon Press, Oxford, second edition, 1986.
\bibitem{Wu} X. Wu, `Distince zeros and simple zeros of Dirichlet $L$-function,' \textit{available at: http://arxiv.org/abs/1206.1679.}

\end{thebibliography}
\end{document}